\documentclass[a4paper,10pt]{article}

\usepackage{amssymb,amsmath,graphics}

\newenvironment{dwd}{\par\noindent{\bf Proof.}}{\par\rightline{$\blacksquare$}}

\newtheorem{theo}{Theorem}
\newtheorem{prop}{Proposition}  
\newtheorem{coro}{Corollary}

\newtheorem{lema}{Lemma}
\newtheorem{defi}{Definition}
\def\be#1\ee{\begin{equation}#1\end{equation}}
\newcommand{\ba}{\begin{eqnarray} }
\newcommand{\ea}{\end{eqnarray} }
\def\bt#1\et{\begin{theo}#1\end{theo}}
\def\bl#1\el{\begin{lema}#1\end{lema}}
\def\bp#1\ep{\begin{prop}#1\end{prop}}
\def\bd#1\ed{\begin{defi}#1\end{defi}}

\def\ccC{{\cal C}}

\def\ccE{{\cal E}}

\def\va{\varepsilon}
\def\ra{\rightarrow}

\def\E{\mathbf{E}}
\def\P{\mathbf{P}}

\def\R{{\mathbb R}}

\def\ls{\leqslant}
\def\gs{\geqslant}
\def\for{\mbox{for}}

\begin{document}

\title{\bf Some remarks on the Sudakov minoration
}
\author{Witold Bednorz
\footnote{Support: Polish Ministry of Science and Higher Education Iuventus Plus Grant no. IP 2011 000171}
\footnote{Institute of Mathematics, University of Warsaw, Banacha 2, 02-097 Warszawa, Poland}}

\maketitle

\begin{abstract}
In this paper we discuss Sudakov type minoration for the dependent setting.
Sudakov minoration is a well known property first proved for centered Gaussian  
processes which states that for well separated points there is a natural lower bound 
on the expectation of the supremum of such a process. We generalize this concept
for the dependent setting where we consider log concave random variables and then 
discuss methods of proving the property.
\end{abstract}

\section{Introduction}\label{sect0}

Consider a random vector $X=(X_1,X_2,...,X_n)$ in $\R^n$ which has log-concave distribution $\mu_X$.
It means that for any non empty measurable sets $A$ and $B$
\be
\mu_X(\lambda A+(1-\lambda)B)\gs \mu_X(A)^{\lambda}\mu_X(B)^{1-\lambda},\;\;\for\;0<\lambda<1.
\ee
Due to the Borel's result \cite{Bor} it means that $\mu$ is supported on the 
affine subspace of $\R^n$ and there exists a density of $\mu$ on the subspace of the form
$\exp(-U(x))$, where $U$ is a convex function. 
\smallskip

\noindent
Consider a finite set $T\subset \R^n$ and a process $(X_t)_{t\in T}$ given by $X_t=\langle t, X\rangle$.
One of the main questions for the analysis of $X$ is to understand the quantity $\sup_{t\in T}X_t$ for
arbitrary set $T\subset \R^n$. The concentration type inequalities describe how well $\sup_{t\in T}X_t$ concentrates
around its mean i.e. $\E\sup_{t\in T}X_t$. In this paper we treat the question what can be said about $\E \sup_{t\in T}X_t$.
\smallskip

\noindent  
We first recall a trivial upper bound on $\E \sup_{t\in T}X_t$. 
\begin{prop}\label{prop0}
Suppose that $X$ satisfies $\|X_t\|_p\ls A$ and $|T|\ls \exp(p)$. 
Then
$$
\E \max_{t\in T}X_t \ls eA
$$
\end{prop}
\begin{dwd}
Obviously
\begin{align*}
&\E\max_{t\in T}X_t \ls \E (\sum_{t\in T}|X_t|^p)^{\frac{1}{p}}\ls (\E \sum_{t\in T} |X_t|^p)^{\frac{1}{p}}\ls\\
&\ls |T|^{\frac{1}{p}}A=eA.
\end{align*}  
\end{dwd} 
The aim of this paper is to reverse the inequality. Obviously it is not possible without additional assumptions on the points in $T$ so
we require that any different points $s,t$ in $T$ are well separated.
The lower bound on $\E \sup_{t\in T}X_t$ under the increment condition is called Sudakov type minoration named after the first result in this direction \cite{Sud} obtained for a vector $X$ of independent centered Gaussian random variables.
Sudakov type minoration is known for independent log concave
random variables and few cases of the general log-concave setting. 
\smallskip

\noindent
We formulate the main problem for this paper. Suppose that $T$ is a finite set. In the optimal case we require that $|T|=\exp(p)$, more generally
we require that there exists a convex increasing function $f:\R_{+}\ra \R_{+}$ such that $f(0)=0$ and $|T|\gs  \exp(f(p))$.
Suppose that we can control all the increments in the following sense
\be\label{n-2}
\|X_t-X_s\|_p\gs A,\;\;\mbox{for all}\;s,t\in T,s\neq t,
\ee
where $p\gs 1$. Note that we can always assume that $p\in Z_{+}$, $p\gs 1$.
In the Sudakov minoration we aim to show that (\ref{n-2}) implies that  
\be\label{n-3}
\E \sup_{t\in T}X_t\gs K^{-1}A,
\ee
where $K$ is an absolute constant. We recall some well known examples when this scheme works:
\begin{enumerate}
\item Gaussian case. Let $X_i=g_i$, $i=\{1,2,...,n\}$ where $g_i$ are independent standard normal variables.
In this case we can apply that $\|X_t-X_s\|_p\sim p^{\frac{1}{2}}\|t-s\|_2$. The meaning of (\ref{n-2}) is that $\|t-s\|_2\gs C^{-1}p^{-\frac{1}{2}}A$
for some absolute constant $C$. Hence by the usual Sudakov minoration (e.g. Theorem 3.18 in \cite{Le-Ta})
$$
\E \sup_{t\in T}X_t \gs K^{-1}p^{\frac{1}{2}} p^{-\frac{1}{2}}A=K^{-1}A. 
$$
\item Bernoulli case. Let $X_i=\va_i$, $i\in \{1,2,...,n\}$, where $\va_i$ are independent random signs, i.e. $\P(\va_i=\pm 1)=\frac{1}{2}$.
Let 
$$
d_p(s,t)=\| \sum^n_{i=1}(t_i-s_i)\va_i\|_p,\;\;\mbox{for}\;s,t\in \R^n.
$$ 
There is a Hitchenko characterization
of $d_p(s,t)$ namely $d_p(s,t)$ is comparable with $(\sum^{p}_{i=1}|t^{\ast}_i|+\sqrt{p}(\sum_{i>p}|t_i^{\ast}|^2)^{\frac{1}{2}})$, where
$|t^{\ast}_1|\gs |t^{\ast}_2|\gs ...\gs |t^{\ast}_n|$ is non increasing rearrangement of $t_1,...,t_n$. For our purposes we need that 
for some absolute constant $C_0\gs 1$
\be\label{n5}
d_{p}(t,0)=\|\sum^n_{i=1}t_i\va_i\|_{p}\gs C_0^{-1}(\sum^{p}_{i=1}|t^{\ast}_i|+p^{\frac{1}{2}}(\sum_{i>p}|t_i^{\ast}|^2)^{\frac{1}{2}}).
\ee
Therefore $d_p(s,t)\gs A$ means that $s-t\not\in AC^{-1}(B^n_1+p^{-\frac{1}{2}}B^n_2)$, where $B^n_q=\{x\in \R^n:\;\sum^n_{i=1}|x_i|^q\ls 1\}$ for any $q\gs 1$.
Consequently either $\E \sup_{t\in T}X_t\gs K^{-1}A$ or $\E \sup_{t\in T}X_t \ls K^{-1}A$ and then by Theorem 4.15 in \cite{Le-Ta}
$$
\E\sup_{t\in T}X_t\gs p^{\frac{1}{2}}(C^{-1}p^{-\frac{1}{2}}A)\gs C^{-1}A.
$$
It implies that $\E \sup_{t\in T} X_t\gs \min\{C^{-1},K^{-1}\}p$.
\item Independent Exponentials. Let $X_i=\xi_i$, where $\xi_i$ are independent symmetric andsuch that 
$\P(|\xi_i|\gs t)=\exp(-C_{\alpha}t^{\alpha})$, where $\alpha\gs 1$. In this setting
the Sudakov type minoration was proved by Talagrand in \cite{Tal1}.
\item Canonical processes. Let $X_i=\xi_i$, where $\xi_i$ are independent symmetric and such that $-\log(\P(|\xi_i|\gs t))$ is convex.
Sudkov minoration for such canonical processes is due Latala \cite{Lat} and is based on a tricky induction argument. 
\end{enumerate}
In this paper we show some progress concerning the dependent case.  We do assume that $X$ is one unconditional. Note that this assumption
helps to reduce the question to a quite natural setting. There are results \cite{La-Pr} that explores the question of general log concave random variables
but they are usually much weaker than what can be proved for the one unconditional case.
\smallskip

\noindent The plan of the paper is the following. In the section \ref{sect1} we collect basic properties of log concave random variables
we need to establish our results. We discuss slightly more general properties than log concavity like $\alpha$ concentration. Section \ref{sect2} concerns 
the main simplification argument which helps to reduce the general question to sets
$T\subset \R^n$ with certain structure of points. In the section \ref{sect3} we explore the case of independent random variables where we recall how the proof 
of the Sudakov minoration works as well as a lot of notation we will need later on. Section \ref{sect4} concerns the case of disjoint supports
where we assume a trivial structure of the simplified set $T$. Finally  in the section \ref{sect5} we study our new setting where the Sudakov type inequality can be proved which is called the common witness existence. 
\smallskip

\noindent Since there will be a lot of constants in the paper we describe our strategy to name them. We use $\alpha,\beta,\gamma,c,C,D,K$ for main constants or constants in the formulation of the results we prove. We use $C_i$, $i=0,1,2,...$ for constants in the assumptions or characterizations which are of meaning for the paper. Finally we use the notation $\bar{C}_i$, $i=0,1,2,...$ for constants inside the proofs we give. Note that for different proofs these constant may vary.

\section{Basic tools}\label{sect1}

We do assume that $X$ is isotropic which mean that $\E X_i=0$ for $i\in\{1,2,...,n\}$ and $\E X_i X_j=\delta_{i,j}$ for $i,j\in \{1,2,...,n\}$. 
In particular it implies that $\E |X_t|^2=\|t\|^2_2$ for any $t\in T$.

\subsection{Bobkov-Nazarov domination}

The first property of log concave $X$ is the Bobkov-Nazarov inequality \cite{Bob}. 
Let $\ccE=(\ccE_1,...,\ccE_n)$ be a vector of independent symmetric exponentials, $\P(|\ccE_i|\gs u)=e^{-u}$, $u\gs 0$.
Bobkov-Nazarov inequality states that tails of $(|X_1|,...,|X_n|)$ are dominated by tails of $(|\ccE_1|,...,|\ccE_n|)$ namely
\begin{eqnarray}\label{bob}
&& \P(\bigcap^n_{i=1} \{|X_i|\gs \sqrt{6} u_i\}) \ls \exp(-\sum^n_{i=1} u_i)=\nonumber\\
&&=\prod^n_{i=1}|P(|\ccE_i|\gs u_i)=\P(\bigcap^n_{i=1}\{|\ccE_i|\gs u_i\}) ,\;\;\mbox{for}\;u_i\gs 0.
\end{eqnarray}
The result (\ref{bob}) is crucial to establish main simplifications of the set $T$. We can slightly 
relax the requirements to enable more general distributions than log-concave which we
need for the independent entries case.
\smallskip

\noindent
We assume that for a given constants $C_1,C_2>0$ and any $I\subset \{1,2,...,n\}$ the following inequality holds
\begin{eqnarray}
&& \P(\bigcap_{i\in I} \{|X_i|\gs C_1 u_i\} )\ls \exp(-\sum_{i\in I} u_i)=\nonumber \\
\label{n2} && =\prod_{i\in I}\P(|\ccE_i|\gs u_i)=\P(\bigcap^n_{i=1}\{|\ccE_i|\gs u_i\}),\;\;\mbox{for}\;u_i\gs C_2.
\end{eqnarray}
For log concave vectors (\ref{n2}) is satisfied with $C_1=\sqrt{6}$ and $C_2=0$.
Note that the condition is slightly less restrictive than the log-concavity. It basically
states that $X$ has sub-exponential distribution for each of its marginals starting from large enough 
arguments. For our purposes we need the following consequence of the property (\ref{n2}) 
\be\label{n3}
\|\sum^n_{i=1}t_i X_i1_{|X_i|\gs C_1C_2}\|_p\ls C_1 \|\sum^n_{i=1} t_i \ccE_i\|_p\;\;\mbox{for all}\;t\in T.
\ee
To get (\ref{n3}) it suffices to show that for any $I\subset \{1,2,...,n\}$ and integers $k_i$, $i\in I$
$$
\E \prod_{i\in I} |X_i|^{2k_i}\prod_{i\in I} 1_{|X_i|>C_1C_2}\ls C_1^{2\sum_{i\in I}k_i}\E \prod_{i\in I} |\ccE_i|^{2k_i}.
$$ 
This is due to (\ref{n2}) and the integration by parts. Recall that the consequence of Kwapien-Gluskin \cite{Kwa} characterization of moments of sums of independent random variables  applied to $(\ccE_1,...,\ccE_n)$ is that 
$$
\|\sum^n_{i=1} t_i \ccE_i\|_p\ls C_3(p\|t\|_{\infty}+\sqrt{p}\|t\|_2),
$$
where $C_3$ is an absolute constant. Therefore using (\ref{n3}) we obtain
\be\label{n4}
\|\sum^n_{i=1}t_i X_i 1_{|X_i|>C_1C_2}\|_p\ls C_1C_3(p\|t\|_{\infty}+\sqrt{p}\|t\|_2).
\ee
The next property we need is the so called $\alpha$ concentration.

\subsection{alfa-Concentration}

The concept is of independent interest and therefore we specify the property for any random vector $X$. We say that $X$ satisfies $\alpha$-concentration, if
there exists a universal constant $\alpha\gs 1$ such that for any norm $\|\cdot\|$ on $\R^n$ there holds
\be\label{n1}
\|X\|_p\ls \alpha\frac{p}{q}\|X\|_q,\;\;2\ls q\ls p,
\ee
where $\|X\|_p=(\E \|X\|^p)^{\frac{1}{p}}$.  For the sake of simplicity we need 
comparison with the first moment of $\|Y\|$. Let
$$
\|X\|_{\varphi_1}=\inf\{C>0:\E \exp(C^{-1}\|X\|)\ls 2\}.
$$
In particular the inequality (\ref{n1}) implies that $\|X\|_{\varphi_1}\ls \alpha C\|X\|_2$, 
for some universal constant $C$, which is
equivalent to
$$
\P(\|X\|\gs  \alpha C u \|X\|_2)\ls e^{-u},\;\;u\gs 0.
$$
The consequence of the $\alpha$ concentration
is the basic control of the distribution of $X_t$, i.e. we have the upper and lower bound
on the tail probability of $\|X\|$:
\begin{enumerate}
\item  Upper tail bound: for any $r\gs 2$ and $u\gs 1$
\be\label{r1}
\P( \|X\|\gs e\alpha u \|X\|_r)\ls \P(\|X\|\gs e\|X\|_{ru})\ls e^{-ru}. 
\ee
\item Lower tail bound: for $r\gs 2$ and $u\in [0,1]$ then for some constant $C\gs 1$ that depends on $\alpha$ only
\be\label{r2}
C \P(\|X\|\gs C^{-1} u\|X\|_r) \gs C\P(\|X\|\gs \alpha C^{-1}\|X\|_{ru})\gs  e^{-ru}.
\ee
\end{enumerate}
All log-concave vectors $X$ in $\R^n$ satisfies this type of concentration with some absolute value of $\alpha$. 
It implies that for a given $\alpha$ variables $X_t=\langle t,X \rangle$, $t\in T$ satisfies 
\be\label{n111}
\| X_t\|_p\ls \alpha \frac{p}{q} \|X_t\|_q,\;\;2\ls q\ls p,
\ee
where $\|X_t\|_p=(\E |X_t|^p)^{\frac{1}{p}}$. Therefore by (\ref{r1}) and (\ref{r2})
we have a control on the tail probability of each $|X_t|$. The slightly more 
involved analysis leads to full understanding of moments.

\subsection{Characterization of moments}

The simplifications we describe in the next section will enable us to consider sets $T$ that contains only
points of thin and different supports. Towards this aim let us introduce the following notation.
For any $t\in \R^n$ we define its support $I(t)\subset \{1,2,...,n\}$ by 
$$
I(t)=\{i\in \{1,2,...,n\}:\;|t_i|>0\}.
$$
Then for any set $J\subset \{1,2,...,n\}$ let us define $t 1_J=(t_i 1_{i\in J})^n_{i=1}$ and
$X_t 1_J=\sum_{i\in J} t_i X_i$ for $t\in \R^n$. Moreover let 
$$
|X|_t=\sum^n_{i=1} |t_i||X_i|\;\;\mbox{and}\;\;|X|_t1_J=\sum_{i\in J}|t_i||X_i|\;\;\mbox{for}\;\;t\in \R^n
$$
be a positive bound on $X_t$.
\smallskip

\noindent
The thin support means that at least $|I(t)|\ls p$. Our basic basic simplification will show that
we can always require that points in $T$ satisfy this requirement. Therefore
to characterize $\|X_t-X_s\|_p=\|X_{t-s}\|_p$ it suffices to bother only the case when $|I(t-s)|\ls p$.
In this setting the following result of Latala \cite{La-Mo} applies.
\bt\label{bambam}
Suppose that $|\mathrm{supp}(t)|\ls p$ then
$$
\|\sum_{i\in I(t)} t_i X_i\|_p\sim \| \sum_{i\in I(t)}t_i|X_i|\|_p \sim\sup\{\sum_{i\in I(t)}|t_i|a_i:\;\P(\bigcap_{i\in I(t)} \{|X_i|\gs a_i\})\gs e^{-p}\}.
$$
It means that there exists $a\in \R^n$ such that $\P(\bigcap_{i\in I(t)} \{|X_i|\gs a_i\})\gs e^{-p}$, $a_i\gs C^{-1}>0$ and
$$
D^{-1}\sum_{i\in I(t)}|t_i|a_i\ls \|\sum_{i\in I(t)} t_i X_i\|_p\ls D \sum_{i\in I(t)}|t_i|a_i, 
$$
where $C,D$  are absolute constants.
\et
It means that to understand $\|\sum_{i\in I(t)} t_iX_i\|_p$ it suffices to consider a witness $a\in \R^n$
supported on $I(t)$ such that $\P(\bigcap_{i\in I(t)} \{|X_i|\gs a_i\})\gs e^{-p}$ that certifies
the linear form to be large, i.e. $\sum_{i\in I(t)} t_i a_i\sim \|X_t\|_p$. We need a
slight improvement of the result. 
\begin{theo}\label{bambam1}
Let $F:\R^n \ra \R_{+}$ be a Borel measurable that satisfies:
\begin{enumerate}
\item $F(0)=0$ and $F(x_1,...,x_n)\ls F(|x_1|,...,|x_n|)$;
\item $F$ restricted to $\R^n_{+}$ is increasing on each coordinate, i.e. for $x\in \R^n_{+}$, $i\in \{1,2,...,n\}$ 
and $\va\gs 0$ there holds $F(x+\va e^i)\gs F(x)$ where $e^i_j=\delta_{i,j}$;
\item in each direction of $\R^n_{+}$ function $F$ satisfies $\triangle_2$ condition in $0$, i.e. there exists $\bar{\alpha}>0$
and $\bar{C}>0$ such that for any $t\in [0,1]$ and $x\in \R^n_{+}$ there holds $F(tx)\gs \bar{C}^{-1}t^{\bar{\alpha}}F(x)$.  
\end{enumerate}
Then for any $p\gs n$ 
$$
\|F(X_1,...,X_n)\|_{p}\sim \sup\{F(a_1,...,a_n):\; \P(\bigcap^n_{i=1}\{|X_i|\gs a_i\})\gs e^{-p}\}.
$$
In particular it means that there exists $a\in \R^n$ such that 
$$
\P(\bigcap^n_{i=1}\{|X_i|\gs a_i\})\gs e^{-p},\;\;\mbox{and}\;\; a_i\gs C^{-1}>0 
$$
and
$$
D^{-1} F(a_1,...,a_n) \ls  \|F(X_1,...,F_n)\|_{p}\ls DF(a_1,...,a_n),
$$
where $C,D$ are universal constants.
\end{theo}
\begin{dwd}
The lower bound is easy. Suppose that there exists $a\in \R^n$ such that
$$
\P(\bigcap^n_{i=1}\{|X_i|\gs a_i\})\gs e^{-p}.
$$
By the one unconditional of $X$
$$
\E |F(X_1,...,X_n)|^p\gs 2^{-n}\E |F(|X_1|,...,|X_n|)|^p\gs (F(a_1,...,a_n))^p 2^{-n}e^{-p}. 
$$
Since $n\ls p$ it implies that
$$
\|F(X_1,...,X_n)\|_{p}\gs 2^{-1}e^{-1} F(a_1,...,a_n).
$$
Therefore
$$
\|F(X_1,...,X_n)\|_{p}\gs 2^{-1}e^{-1}\sup\{F(a_1,...,a_n):\; \P(\bigcap^n_{i=1}\{|X_i|\gs a_i\})\gs e^{-p}\}.
$$
To prove the upper bound we need the main tool of \cite{La-Mo}. W.l.o.g. we may assume that there exists 
a non-degenerate density $e^{-U(x)}$ of $X$. Let
$$
K_p=\{y\in \R^n:\; U(y)-U(0)\ls p\},\;\;\|F\|_{K_p^{\circ}}=\sup\{F(y_1,...,y_n):\;y\in K^p\}
$$ 
It is proved in \cite{La-Mo} that there exists an absolute constant $\bar{C}_0$ such that
$$
\P(X\in \bar{C}_0 K_p)\gs 1-e^{-p}.
$$
Clearly
$$
\P(F(X_1,...,X_n)> \bar{C}_0 \|F\|_{K^{\circ}_p}) \ls \P(X\not\in \bar{C}_0 K_p)\ls e^{-p}
$$
Moreover by the log-concavity of $X$ and $t\gs 1$
$$
\P(F(X_1,...,X_n)> \bar{C}_0 t\|F\|_{K^{\circ}_p})\ls \P(X\not\in t\bar{C}_0 K_p)\ls e^{-tp}.
$$
Integration by parts implies that
$$
\|F(X_1,...,X_n)\|_p\gs D^{-1}\|F\|_{K_p^{\circ}},
$$
where $D$ is a universal constant. It suffices to choose $y\in \R^n$ such that $U(y)-U(0)\ls p$ and
$F(y)\gs D^{-1} \|F\|_{K_p^{\circ}}$. We finish the proof in the same way as in Corollary 2 from \cite{La-Mo}.
First is is easy to notice that $U(0)\ls \frac{3}{2}p$ and hence $U(y)\ls \frac{5p}{2}$.
Then the basic properties of log-concave vectors imply that for a universal constant $\bar{C}_1$ 
$U(\bar{C}^{-1}_1,...,\bar{C}^{-1}_1)\ls \frac{5}{2}p$. Hence for $z_i=\frac{1}{2}(\bar{C}^{-1}_1+y_i)$, $i\in \{1,2,...,n\}$
we have that $U(z)\ls \frac{5}{2}p$ and consequently using that $U$ is coordiante increasing
$$
\P(\bigcap^n_{i=1}\{X_i\gs \frac{z_i}{2}\})\gs e^{-U(z)}\prod^n_{i=1}\frac{z_i}{2}\gs e^{-\frac{5p}{2}}(4\bar{C}_1)^{-n}.
$$
Since $p\gs n$ and $s\ra -\ln \P(X_1\gs s_1,...,X_n\gs s_n)$ is convex we get
$$
\P(\bigcap^n_{i=1}\{|X_i|\gs \bar{C}^{-1}_2 z_i\})=2^n\P(\bigcap^n_{i=1}\{X_i\gs \bar{C}^{-1}_2z_i\})\gs e^{-p}
$$
for sufficiently large $\bar{C}_2$. By the properties of $F$
$$
F(\bar{C}^{-1}_2 z_1,....,\bar{C}^{-1}_2z_n)\gs \bar{C}^{-1} \bar{C}_2^{-\bar{\alpha}} F(y_1,...,y_n)\gs D^{-1}\bar{C}^{-1}\bar{C}^{-\bar{\alpha}}_2 p.
$$
which ends the proof.
\end{dwd}
\begin{coro}\label{common}
Suppose that $|I(t)|\ls p$ then for any class $\ccC$ of subsets of $I(t)$ the following holds
$$
\|\min_{C\in \ccC}|\sum_{i\in C} t_i X_i|\|_p\sim \sup\{\min_{C\in \ccC}\sum_{i\in I(t)} |t_i|a_i:\;\P(\bigcap_{i\in I(t)}\{|X_i|\gs a_i\})\gs e^{-p}\}.
$$
Consequently there exists $a\in \R^n$ supported in $I(t)$ such that 
$$
\P(\bigcap_{i\in I(t)} \{|X_i|\gs a_i\})\gs e^{-p}
$$
and
$$
D^{-1}\|\min_{C\in \ccC}|\sum_{i\in C} t_i X_i|\|_p\ls \min_{C\in \ccC}|\sum_{i\in C}|t_i|a_i|\ls D \|\min_{C\in \ccC}|\sum_{i\in C} t_i X_i|\|_p
$$
\end{coro}
\begin{dwd}
First note that due to one unconditionality of $X$ we may assume that $t_i\gs 0$ for $i\in I(t)$. 
Then it suffices to define $F(x)=\min_{C\in \ccC}|\sum_{i\in C}t_i x_i|$ and use Theorem \ref{bambam1}.
\end{dwd} 
In this way we obtain the tool for the so called common witness existence. The point is that
if we have that for a class $\ccC$ we can show that $\|\min_{C\in \ccC}|\sum_{i\in C} t_i X_i|\|_p$
is greater then $A$ then we have a witness $a\in \R^n$ which is good for any subset $C\in \ccC$.
We apply the result to $\ccC(t)=\{I(t)\backslash I(s):\; \|X_t 1_{I(t)\backslash I(s)}\|_p\gs A\}$
for each $t\in T$.

\subsection{Exponential inequality}

One of the most powerful tool for log-concave random variables are exponential type inequalities.
Let $X=(X_1,...,X_n)$ be log-concave. We say that $X$ satisfies exponential concentration
with constant $\beta$, i.e. whenever $\P(X\in B)\gs \frac{1}{2}$ for a Borel set $B$ then
\be\label{n6}
\P(X\in B+\beta u B^n_2)\gs 1-e^{-u},\;\;\mbox{for}\;u>0,
\ee
where $B^n_2=\{x\in \R^n:\; \sum^n_{i=1}x_i^2\ls 1 \}$.
For log concave vectors this inequality holds at least with $\beta\ls Cn^{\frac{1}{2}-\va}$ for some $\va>0$, e.g. $\va =\frac{1}{8}$. In the next section we will need the optimal known estimate \cite{Kla} for $\beta$ under the one unconditionality assumption, i.e. $\beta\ls C\log n$. In general it is conjectured that (\ref{n6}) holds with $\beta$ which does not depend on $n$ -KLS conjecture \cite{KLS}.
\smallskip

\noindent
The exponential inequality gives some geometrical understanding of the distribution of $X$. We use the idea
to first give a new proof of the Sudakov minoration for disjoint supports i.e. when $I(t)\cap I(s)=\emptyset$
for all $s,t\in T$ and $s\neq t$. Then we show that the argument can be slightly generalized to the case when 
the common witness exists for each $t\in T$.  We conclude that the Sudakov minoration for $T$ holds at least when $f(p)=p^2$
and sometimes this can be improved to $f(p)=p\log(1+p)$.

\section{How to simplify the problem}\label{sect2}

Assume that $X=(X_1,...,X_n)$ is isotropic and one-unconditional. In this section we analyze a list of simplifications
of the setting in which Sudakov minoration has to be proved. Recall that although the best form of the 
Sudakov minoration works for $|T|\gs \exp(p)$ we consider much more general requirement that $|T|\gs \exp(f(p))$, 
where $f$ is increasing and $f(0)=0$.
\smallskip

\noindent
Our first observation is that one can always require that $0\in T$. This due to isotropy, i.e. for any $s\in T$
we have
$$
\E \sup_{t\in T}X_t=\E \sup_{t\in T}X_t-X_s=\E \sup_{t\in T}X_{t-s}=\E \sup_{t\in T-s}X_t.
$$ 
By the symmetry of $X_t$ it implies that
\be\label{mamba}
\E \sup_{t\in T} X_t=\E \sup_{t\in T}(X_t)_{+}\gs \frac{1}{2}\E \sup_{t\in T}|X_t|.
\ee
Therefore to get $\E \sup_{t\in T} X_t\gs A$ it suffices to prove $\E \sup_{t\in T}|X_t|\gs A$.
Due to the homogeneity of the problem we may require that $A=p$, which means
that (\ref{n-2}) can be rewritten as
$$
\|X_t-X_s\|_p\gs p,\;\;\mbox{for all}\;s,t\in T,s\neq t
$$
and (\ref{n-3}) in the view of (\ref{mamba}) as
$$
\E \sup_{t\in T}|X_t|\gs K^{-1}p,
$$
where $K$ is a universal constant that depends on the function $f$ only. 
\smallskip

\noindent
We are ready to present more involved simplifications of the set $T$. Towards this
aim we have to assume some regularity of the distribution of $X$.
The fact that $X$ is one unconditional implies that we can benefit from 
the tools invented for Bernoulli random variables (see chapter 4 in \cite{Le-Ta}). 
On the other hand we need a control from above on tails of $X_t-X_s$, for $s,t\in T$.
As we have mentioned in the previous section tails of log-concave vectors are dominated
by independent symmetric exponentials. For our purposes we need a slightly weaker form
of this property i.e. we assume (\ref{n2}) which implies its useful consequences (\ref{n3}) and (\ref{n4}).
\smallskip

\noindent
We prove two results. The first one concerns the perfect case of $f(p)=p$.
The important feature of the proof is that it indeed requires the exponential number of points in $T$.
\begin{prop}\label{prop1}
Suppose that $X$ satisfies (\ref{n2}), then for any $T$ such that $0\in T$, $|T|=\exp(p)$ and
$$
\|X_t-X_s\|_p\gs p,\;\;s,t\in T,\;s\neq t
$$
to prove that for a universal $K$ 
\be\label{n5-1}
\E \sup_{t\in T}|X_t|\gs K^{-1} p
\ee
it suffices to show that for a suitably small $\delta$ and $p$ suitably large 
there exists a universal constant $K$ such that
for any set $T$ that satisfies:
\begin{enumerate}
\item  $|T|\gs \exp(\frac{1}{4}p)$ and $0\in T$;
\item  for each $i\in \{1,2,...,n\}$
$$
t_i\in \{0,k_i\},\;\;\mbox{where}\;k_i\gs \rho;
$$
\item  for each $t\in T$ 
\be\label{n50}
\sum^n_{i=1} k_i 1_{t_i\neq 0}\ls 2C_0\delta p,
\ee
where $\rho\ls e^{-1}$ and $\rho/\log \frac{1}{\rho}= 4C_0\delta$ and $C_0$ is from (\ref{n5});
\item  for all $s,t\in T$, $s\neq t$
$$
\|X_t-X_s\|_p\gs \frac{p}{2};
$$ 
\end{enumerate}
the following inequality holds
$$
\E\sup_{t\in T}|X_t|\gs K^{-1}p.
$$
\end{prop}
\begin{dwd}
The proof is based on the number of straightforward simplifications. 
\smallskip

\noindent 
{\bf Step 1} Recall that $d_p(t,s)=\|\sum^n_{i=1} (t_i-s_i)\va_i\|_p$ for $s,t\in T$.
We may assume that $p\gs 1$ is suitably large. Moreover we may consider $T$ such that $0\in T$, $|T|\gs  \exp(\frac{3p}{4})$ and
$d_p(t,0)\ls \delta p$ for all $t\in T$, where $\delta\ls 1$ can be suitably small.
\smallskip

\noindent
Obviously it suffices to prove the result for $p$ that are sufficiently large. Let $N(T,d_p,u)$ is the entropy number for $T$ i.e.
the minimal cardinality of balls of radius $u$ in $d_p$ distance that are required to cover $T$.
As we have already mentioned by the Talagrand's \cite{Tal1} result (e.g. Theorem 4.15 in \cite{Le-Ta}) if $N(T,d_p,u)\gs \exp(\frac{p}{4})$ 
then $\E \sup_{s,t\in T} \sum^n_{i=1}(t_i-s_i)\va_i\gs K^{-1}u$, for
a universal $K$. Thus we may assume that $N(T,d_p,\frac{1}{2}\delta p)\ls \exp(\frac{p}{4})$. It implies that
there exists $t_0\in T$ such that 
$$
\{t\in T:\; d_{p}(t,t_0)\ls \delta p\}|\gs |T|\exp(-\frac{p}{4})\gs \exp(-\frac{3p}{4}).
$$
Therefore we may consider set $T'=\{t-t_0: d_{p}(t,t_0)\ls \delta p\}$, which satisfies
all the requirements. 
\smallskip

\noindent
{\bf Step 2} 
Let $\rho\ls 4C_0 \delta\ls e^{-1}$.
We may assume that $0\in T$, $|T|\gs \exp(\frac{p}{4})$ and additionally  
$$
t_i\in (k_i-\rho,k_i+\rho)\cup (-\rho,\rho )\;\;\mbox{for all}\;t\in T\;\;\mbox{and}\;\; 1\ls i\ls n,
$$ 
where $k_i$ are given numbers such that $k_i\gs \rho$ and $\rho\ls e^{-1}$ that satisfies 
$$
\rho/\log \frac{1}{\rho}=4C_0 \delta\ls e^{-1},
$$ 
and $C_0$ is the constant in (\ref{n6}). 
\smallskip

\noindent
Indeed consider measure  $\mu=\otimes^n_{i=1}\mu_i$, where $\mu_i(dx)=\frac{1}{2}e^{-|x|}dx$ for
all $1\ls i\ls n$. For any $x\in \R^n$ and $t\in T$
$$
T_x=\{t\in T:\;t_i\in (x_i-\rho,x_i+\rho)\cup (-\rho,\rho),\;i=1,2,...,n\}
$$
and
$$
A_t=\{x\in \R^n:\;t_i\in (x_i-\rho,x_i+\rho)\cup (-\rho,\rho),\;i=1,...,n\}.
$$
Now there are two possibilities either
$$
\mu_i(\{x_i:\;t_i\in (x_i-\rho,x_i+\rho)\cup (-\rho,\rho) \})\gs \rho e^{-|t_i|-\rho}
$$  
or
$$
\mu_i(\{x_i:\;t_i\in (x_i-\rho ,x_i+\rho)\cup (-\rho,\rho) \})=1,\;\;\mbox{if}\;|t_i|<\rho.
$$
By (\ref{n5}) we get
$$
d_{p}(t,0)=\|\sum^n_{i=1}t_i\va_i\|_{p}\gs C_0^{-1}(\sum^{p}_{i=1}|t^{\ast}_i|+\sqrt{p}(\sum_{i>p}|t_i^{\ast}|^2))^{\frac{1}{2}}.
$$
Therefore since $d_{p}(t,0)\ls \delta p$ and $\rho/\log \frac{1}{\rho}=4C_0\delta$, $\rho\ls e^{-1}$ we get 
\be\label{n56}
|\{i\in \{1,...,n\}:\;|t_i|\gs \rho\}|\ls \frac{p}{4\log \frac{1}{\rho}}\ls \frac{p}{4}.
\ee
It implies also 
\be\label{rose}
\sum^n_{i=1} |t_i|1_{|t_i|\gs \rho}\ls C_0\delta p\ls \frac{p}{4e}.
\ee
Consequently
$$
\mu(A_t)\gs \rho^{\frac{p}{4\log \frac{1}{\rho}}} \exp(-\frac{p}{4e})\exp(-\rho\frac{p}{4} )\gs \exp(-\frac{p}{2}).
$$
However using that $|T|\gs \exp(\frac{3p}{4})$ we obtain
$$
\int \sum_{t\in T}1_{A_t}(x)\mu(dx)\gs |T|\exp(-\frac{p}{2} )\gs \exp(\frac{p}{4}).
$$
Therefore we get that there exists
at lest one point $k\in \R^n$ such that
$$
|T_k|\gs \exp(\frac{p}{4}).
$$  
It is obvious that $|k_i|$ may be chosen in a way that $|k_i|\gs \rho$.
Using (\ref{rose}) together with $|k_i|\gs \rho$ and $t_i\in (k_i-\rho,k_i+\rho)$ we obtain
$$
\sum^n_{i=1} (|k_i|-\rho)\vee \rho 1_{|t_i|\gs \rho}\ls C_0\delta p.
$$
Clearly $(|k_i|-\rho)\vee \rho\gs \frac{1}{2}|k_i|$ and therefore
$$
\frac{1}{2} \sum^n_{i=1}|k_i|1_{|t_i|\gs \rho}\ls C_0\delta p,
$$
which implies (\ref{n50}). Clearly by the one unconditionality of $X$ we may consider only $k_i$, $i\in \{1,2,...,n\}$ positive. 
\smallskip

\noindent
{\bf Step 3} It suffices to consider set $T$ which additionally satisfies $t_i\in \{0,k_i\}$ where $k_i\gs \rho$ and still
$$
\|X_t-X_s\|_p\gs \frac{p}{2},\;\;\mbox{for all}\;s,t\in T,s\neq t.
$$
Consider the following function
$$
\varphi_i(t_i)=\left\{\begin{array}{lll}
0 & \mbox{if} & |t_i|<\rho \\
k_i & \mbox{if} & |t_i|\gs \rho 
\end{array}\right.
$$
Let $\varphi(t)=(\varphi_i(t_i))^n_{i=1}$. We show that $\| X_{\varphi(t)}-X_{\varphi(s)}\|_p\gs \frac{p}{2}$. 
It requires  some upper bound on $\|X_{t-\varphi(t)}\|_p$. Consider any $s\in T$ then using (\ref{n2}) (or rather (\ref{n4}))
\begin{align*}
& \| X_s\|_p\ls  C_1C_2d_p(s,0)+ (\E |\sum^n_{i=1} s_i X_i|^p\prod^n_{i=1}1_{|X_i|\gs C_1 C_2})^{\frac{1}{p}}\ls \\
&  \ls C_1 C_2d_p(s,0)+C_1C_3(p\|s\|_{\infty}+\sqrt{p}\|s\|_2).
\end{align*}
Using the contraction principle (e.g. Theorem 4.12 \cite{Le-Ta}) for $s=t-\varphi(t)$
$$
(p\|t-\varphi(t)\|_{\infty}+\sqrt{p}\|t-\varphi(t)\|_2)\ls 2\rho p+d_p(t,0).
$$
Consequently using (\ref{n5})
\begin{align*}
&\|X_{t-\varphi(t)}\|_p 2C_1C_3 \rho p+(C_0C_1C_3+2C_1 C_2)d_p(t,0)\ls\\
&\ls  (2C_1C_3 \rho+(C_0C_1C_3+2C_1C_2)\delta)p\ls \frac{p}{4}.
\end{align*}
for suitably small $\delta$. Therefore
$$
\|X_{\varphi(t)}-X_{\varphi(s)}\|_{p}\gs \|X_t-X_s\|_{p}-\|X_t-X_{\varphi(t)}\|_{p}-\|X_s-X_{\varphi(s)}\|_{p}\gs \frac{p}{2}.
$$
Suppose we can prove the main result for the constructed set $T$ (of cardinality $\exp(\frac{p}{4})$), namely suppose that
$$
\E \sup_{t\in T}|X_{\varphi(t)}|\gs K^{-1}p,
$$  
for some universal $K$. Recall that 
$$
\|X_{t-\varphi(t)}\|_p\ls (2C_1C_3 \rho+(C_0C_1C_3+2C_1C_2)\delta)p
$$
and therefore by Proposition \ref{prop0} we get
$$
\E\sup_{t\in T} |X_{t-\varphi(t)}|\ls e^{\frac{1}{4}}(2C_1C_3 \rho+(C_0C_1C_3+2C_1C_2)\delta)p.
$$ 
Thus
$$
\E \sup_{t\in T} X_t=\E \sup_{t\in T} X_{\varphi(t)}+X_{t-\varphi(t)}\gs \E \sup_{t\in T}X_{\varphi(t)}- \E\sup_{t\in T} X_{t-\varphi(t)}\gs K^{-1}\frac{p}{2},
$$
for suitably small $\delta$, i.e. $2(2C_1C_3 \rho+(C_0C_1C_3+2C_1C_2)\delta)\ls K^{-1}$.
\end{dwd}
\begin{coro}\label{coro1}
Note that in particular after the simplification points in $T$ are of thin and different support, i.e. $|I(t)|\ls cp$, where $c$ is sufficiently small
and $I(t)\neq I(s)$ if $s\neq t$, $s,t\in T$.
\end{coro}
\begin{dwd}
To see the that supports are thin it suffices to use (\ref{n56}) and observe that it implies
$$
|I(t)|\ls \frac{p}{4\log \frac{1}{\rho}}=cp,\;\;\mbox{where}\;c=\frac{1}{4\log \frac{1}{\rho}}.
$$ 
The supports are different since $t_i\in \{0,k_i\}$ for all $i\in \{1,2,...,n\}$ and $\|X_t-X_s\|\gs \frac{p}{2}$ for all $s\neq t$, $s,t\in T$.
\end{dwd}
Unfortunately if we increase the number of  points in $T$ beyond $\exp(p)$ the above simplification is no longer possible and 
what fails is the last step of the proof where we have to bound the subtracted process.
On the other hand without this step still we can obtain a similar result with a little less control on the structure of points in $T$.
\begin{prop}\label{prop2}
Suppose that $X$ satisfies (\ref{n2}), then  for any $T$ such that $0\in T$, $|T|\gs \exp(f(p))$  
and
$$
\|X_t-X_s\|_p\gs p,\;\;s,t\in T,\;s\neq t
$$
the following inequality holds
$$
\E \sup_{t\in T}|X_t|\gs K^{-1} p
$$
if for a suitably small $\delta$ and sufficiently large $p$ there exists a universal constant $K$ such that for 
any set $T$ of properties:
\begin{enumerate}
\item $|T|\gs \exp(-\frac{3}{4}p+f(p))$, $0\in T$;
\item for each $i\in \{1,2,...,n\}$
$$ 
t_i\in ((k_i-2\rho)_{+} ,k_i)\cup \{0\},\;\;\mbox{where}\; k_i\gs \rho ; 
$$
\item for each $t\in T$
\be\label{brum1}
\sum^n_{i=1} k_i 1_{t_i\neq 0}\ls 2C_0\delta p,
\ee
where $\rho\ls e^{-1}$ and $\rho/\log \frac{1}{\rho}= 4C_0\delta$, for some $C_0\gs 1$;
\item for all $s,t\in T$, $s\neq t$ 
$$
\|X_t-X_s\|_p\gs \frac{p}{2};
$$ 
\end{enumerate}
the following inequality holds
$$
\E\sup_{t\in T}|X_t|\gs K^{-1}p.
$$
\end{prop}
\begin{dwd}
We follow the steps in the proof of Proposition \ref{prop1}.
\smallskip

\noindent 
{\bf Step 1} We may assume that $T$ is of the form $0\in T$, $|T|\gs  \exp(-\frac{p}{4}+f(p))$ and
$d_p(t,0)\ls \delta p$ for all $t\in T$, where $\delta\ls 1$ can be suitably small.
The proof is the same as in Proposition \ref{prop1} and is based on the fact that
if $N(T,d_p,\frac{1}{2}\delta p)\gs \exp(\frac{p}{4})$ then $\E \sup_{t\in T}|X_t|\gs K^{-1}p$ for some universal 
constant $K$.
\smallskip

\noindent
{\bf Step 2} 
Let $\rho= 4C_0 \delta\ls e^{-1}$.
We may assume that $0\in T$, $|T|\gs \exp(-\frac{3p}{4}+f(p))$ and additionally  
$$
t_i\in (k_i-\rho,k_i+\rho)\cup (-\rho,\rho )\;\;\mbox{for all}\;t\in T\;\;\mbox{and}\;\; 1\ls i\ls n,
$$ 
where $k_i$ are given numbers such that $k_i\gs \rho$ and $\rho\ls e^{-1}$ that satisfies 
$$
\rho\log \frac{1}{\rho}=4C_0 \delta\ls e^{-1},\;\;\mbox{where}\; C_0\gs 1
$$ 
is a universal constant. In particular it means that if $|t_i|\gs \rho$ then
$t_i>0$. Moreover (\ref{brum1}) holds, i.e.
$$
\sum^n_{i=1} k_i 1_{|t_i|\gs \rho}\ls 2C_0\delta p.
$$
Again the proof is the same as in Proposition \ref{prop1}.
\smallskip

\noindent
{\bf Step 3} It suffices to consider set $T$ such that
$$
t_i\in ((k_i-2\rho)_{+} ,k_i)\cup \{0\},\;\;\mbox{where}\; k_i\gs \rho,
$$
and at least
$$
\|X_t-X_s\|_p\gs \frac{p}{2},\;\;\mbox{for all}\;\;s,t\in T,s\neq t.
$$
Instead of the function $\varphi$ as in the proof of Proposition \ref{prop1} 
we can use the following 
$$
\varphi_i(t_i)=\mathrm{sign}(t_i)(|t_i|-\rho)_{+}.
$$
Then the Bernoulli comparison (see Theorem 4.15 in \cite{Le-Ta}) follows that 
$$
\E \sup_{t\in T}|X_t| \gs \E \sup_{t\in T}|X_{\varphi(t)}|.
$$
Finally we show 
$$
\| X_{\varphi(t)}-X_{\varphi(s)}\|_p\gs \frac{p}{2}.
$$
Following the proof of Proposition \ref{prop1} we get
for sufficiently small  $\delta$ that
$$
\|X_{t-\varphi(t)}\|_p\ls (2C_1C_3 \rho+(C_0C_1C_3+2C_1C_2)\delta)p\ls \frac{p}{4},
$$
and hence
$$
\|X_{\varphi(t)}-X_{\varphi(s)}\|_{p}\gs \|X_t-X_s\|_{p}-\|X_t-X_{\varphi(t)}\|_{p}-\|X_s-X_{\varphi(s)}\|_{p}\gs \frac{p}{2}.
$$
It completes the proof of the result.
\end{dwd}
\begin{coro}\label{cor2}
After the simplification from Proposition \ref{prop2} points $t\in T$ are of short and different supports.
Namely $I(t)\neq I(s)$ for all $s\neq t$, $s,t\in T$ and $|I(t)|\ls cp$, where $c$ is sufficiently small and absolute constant. 
\end{coro}
\begin{dwd}
We use Proposition \ref{prop2} and hence by (\ref{brum1})
$$
\rho |I(t)|\ls \sum^n_{i=1} k_i 1_{t_i\neq 0}\ls 2C_0\delta p=\frac{p}{2\log \frac{1}{\rho}}.
$$
Therefore $|I(t)|\ls cp$, where $c=1/(2\log \frac{1}{\rho})$.
\smallskip

\noindent
Suppose that after the simplification of Proposition \ref{prop2} we have two points $s,t\in T$, $s\neq t$ of the same 
supports, i.e. $I(s)=I(t)$. By the contraction principle (e.g. Theorem 4.12 \cite{Le-Ta})
\begin{align*}
& \frac{p}{2}\ls \| X_s-X_t\|_p\ls  C_1C_2d_p(s,t)+ (\E |\sum^n_{i=1} (s_i-t_i) X_i1_{|X_i|\gs C_1 C_2}|^p)^{\frac{1}{p}}\ls \\
&  \ls (C_1 C_2)d_p(s,t)+(C_1C_3)(p\|s-t\|_{\infty}+\sqrt{p}\|s-t\|_2)
\end{align*}
and $d_p(s,t)\ls 2\delta p$. By the construction $I(t)=I(s)$ implies $\|t-s\|_{\infty}\ls 2\rho$, moreover
as we have proved $|I(t)|,|I(s)|\ls cp$ and hence 
$$
\|t-s\|_{2}\ls 2\rho\sqrt{cp}.
$$ 
This leads to contradiction for suitably small $\rho$.
\end{dwd}
\begin{coro}\label{cor3}
Suppose that $\delta$ is suitably small. For each $s,t\in T$, $s\neq t$
\be\label{nuj}
 \frac{1}{2}\|\sum_{I(t)\bigtriangleup I(s)}k_i 1_{k_i> 4\rho}X_i\|_p\ls   
 \| X_t -X_s\|_p\ls  2 \| \sum_{I(t)\bigtriangleup I(s)}k_i 1_{k_i> 4\rho} X_i \|_p.
\ee
Moreover for each $s,t\in T$, $t\neq s$ 
\be\label{n15}
2^{-1}C_0^{-1}\|\sum_{I(t)\bigtriangleup I(s)}k_i 1_{k_i> 4\rho} |X_i|\|_p\ls    \| X_t -X_s\|_p\ls  
2\| \sum_{I(t)\bigtriangleup I(s)}k_i 1_{k_i> 4\rho} |X_i| \|_p.
\ee
\end{coro}
\begin{dwd}
To prove the first assertion we observe that by the same argument as in Corollary \ref{cor2}
\begin{eqnarray}
&&\|X_t-\sum_{i\in I(t)} k_i X_i\|_p+\|X_s-\sum_{i\in I(s)} k_i X_i\|_p=
 \|\sum_{i\in I(t)}(t_i-k_i)X_i\|_p+ \nonumber \\
\label{notnow} &&+\|\sum_{i\in I(s)}(s_i-k_i)X_i\|_p\ls  \bar{C}_0^{-1}p\ls 2\bar{C}_0^{-1}\|X_t-X_s\|_p
\end{eqnarray}
where $\bar{C}_0$ is suitably small. In the same we get
$$
\|\sum_{i\in I(t)\bigtriangleup I(s)}k_i 1_{k_i\ls 4\rho}X_i\|_p\ls
\bar{C}_1 p,
$$
where $\bar{C}_1$ can be sufficiently small. It implies (\ref{nuj}) and also the upper bound in (\ref{n15}).
\smallskip

\noindent 
We turn to prove the lower bound in (\ref{n15}).
As we have shown
$$
\|X_t-X_s\|_p\gs 2^{-1}\|\sum_{i\in I(t)\bigtriangleup I(s)}k_i1_{k_i> 4\rho} X_i \|_p.
$$
By (\ref{n5}) we get
$$
\|\sum_{i\in I(t)\bigtriangleup I(s)}k_i1_{k_i>4\rho} X_i \|_p\gs \bar{C}_0^{-1}\|\sum_{I(t)\bigtriangleup I(s)}k_i |X_i|\|_p. 
$$
\end{dwd}
The meaning of the above result is that the simplifications from Proposition \ref{prop1} and \ref{prop2} are
are of the similar power. On the other hand only for $f(p)=p$ we can prove that only supports matters.
In Proposition \ref{prop1} it is obvious. In the general setting of Proposition \ref{prop2} we can 
get the similar result by the following fact. 
\begin{lema}\label{lema5}
Suppose that $f(p)= Cp$ for a given constant $C\gs 1$. Then 
$$
\E \sup_{s,t\in T}|X_t-X_s|1_{I(t)\cap I(s)}\ls D^{-1}p
$$
for sufficiently large $D$.
\end{lema}
 \begin{dwd}
To prove the result first note that by the same argument as in Corollary \ref{cor2} i.e. (\ref{notnow}) we get
$$
\|(X_t-X_s)1_{I(t)\cap I(s)}\|_p \ls \bar{C}_0^{-1}p,
$$
for $\bar{C}_0$ sufficiently large (depending on $\delta$). 
Therefore by Proposition \ref{prop0} (note that $T\times T$ counts not more than $\exp(2 Cp)$ elements) 
$$
\E \sup_{s,t\in T}|X_t-X_s|1_{I(t)\cap I(s)}\ls 2C \bar{C}_0^{-1}p
$$
and hence the result for $D =(2 C)^{-1}\bar{C}_0$.
\end{dwd}
Due to the triangle inequality it implies that whenever it is possible to prove 
$$
\E \sup_{t\in T} \sup_{s\in T}|X_t| 1_{I(t)\backslash I(s)}\gs K^{-1}p
$$
then also $\E \sup_{t\in T}|X_t|\gs 2^{-1}K^{-1}p$ which is
difficult to get without this tool.

\section{Independent entries}\label{sect3}

Let $X$ be isotropic and one unconditional in $\R^n$. 
In this section we assume independence of entries of $X$ as well as $\alpha$ concentration of each of them.
It means we assume that $X=(X_1,...,X_n)$ is such that $X_1,...,X_n$ are independent symmetric and satisfy
$$
\| X_i\|_p \ls \alpha \frac{p}{q} \|X_i\|_q,\;\; p\gs q\gs 1,\;\;i\in\{1,2,...,n\}.
$$
These assumptions enable us to sufficiently control the distribution of $X$ by the independence and (\ref{r1}) and (\ref{r2})
applied to each $X_i$ for $1\ls i\ls n$. Note that the upper tail bound (\ref{r1}) applied for $r=2$ together with the independence of entries 
imply (\ref{n2}) with $C_1=C_2=\alpha e$. It was the main purpose of formulating the weaker form of
the Bobkov-Nazarov tail domination. The consequence of (\ref{n2}) is that we can apply Proposition \ref{prop1}.
\smallskip

\noindent
The meaning of Proposition \ref{prop1} is that we can analyze $T$ such that $|T|\gs \exp(\frac{p}{4})$ which contains only
points of short and different supports, such that $t_i\in \{0,k_i\}$ for $i\in \{1,2,...,n\}$ and
$$
\|X_t-X_s\|_p\gs \frac{p}{2}\;\;\mbox{for all} \;s,t\in T.
$$
For log concave vectors Theorem \ref{bambam} helps to fully characterize $\|X_t-X_s\|_p$.
It basically states that for each $s,t\in T$ there exists a witness $a\in \R^n$ such that 
$\P(\bigcap \{|X_i|\gs a_i\})\gs e^{-p}$ and  
$$
\sum_{i\in I(t)\bigtriangleup I(s)}x_i a_i\gs D^{-1}\frac{p}{2}.
$$
for a universal constant $D$. For the proof of the Sudakov minoration it suffices
to use a global upper bound $r\in \R^n$ for such class of $a\in \R^n$ obtained for all $s,t\in T$. For log-concave vectors we could use the density the density $U(x)=\sum^n_{i=1}U_i(x_i)$ of $\mu_X$ and define $r=(r_i)^n_{i=1}$ as the solution of $U^{\ast}_i(k_i)=k_i r_i-U_i(r_i)$, where
$U^{\ast}_i$ is the conjugate function to $U_i$.
We can slightly generalize this idea using moments of $X_i$, $i\in \{1,2,...,n\}$ which better matches the setting of $\alpha$ concentration.
\smallskip

\noindent 
The point is that such a witness can be to some extent defined by the analyze of single $X_i$, $1\ls i\ls n$.
Fix constant $\gamma>0$. For each $i\in \{1,2,...,n\}$ we define  $r_i\in [2,p]$ by 
\begin{enumerate}
\item $r_i=2$ if $k_i\|X_i\|_2\ls 2\gamma$;
\item otherwise $r_i=p$ if $k_i\|X_i\|_{r}> r\gamma $ for all $r\in [2,p]$  ;
\item otherwise $r_i=\inf\{r\in[1,p]:\;k_i\|X_i\|_{r}=r\gamma\}$.
\end{enumerate}
Obviously one of the above three possibilities must hold. We state the crucial consequence of
the condition $\|X_t-X_s\|_{p}\gs \frac{p}{2}$ for the independent case.
\begin{lema}\label{lema3}
Fix $\gamma=(8\alpha e)^{-1}$. Then for all $s,t\in T$ after the simplification from Proposition \ref{prop1}
the following inequality holds
\be\label{inqu}
\sum_{i\in I(t)\bigtriangleup I(s)}r_i1_{r_i> 2}\gs \gamma p.
\ee
\end{lema}
\begin{dwd}
The easy case is when there exist at least one $i\in I(t)\bigtriangleup I(s) $ such that
$r_i=p$ in which (\ref{inqu}) trivially holds. Thus we can assume that $r_i<p$ for all $i\in I(t)\bigtriangleup I(s)$ which implies
by the construction of $r_i$
\be\label{tlen}
k_i|X_i|= \frac{k_i\|X_i\|_{r_i}}{r_i}\frac{r_i|X_i|}{\|X_i\|_{r_i}}\ls \alpha \gamma e  \frac{r_i|X_i|}{e\alpha\|X_i\|_{r_i}},
\ee
By (\ref{r2}) we have that at least for $t\gs r_i$, random variable $\frac{r_i|X_i|}{\alpha e\|X_i\|_{r_i}}$ has its tail dominated by $|\ccE_i|$.
Therefore  
\begin{align*}
&\|X_t-X_s\|_p\ls (\alpha e \sum_{i\in I(t)\bigtriangleup I(s)} k_i\|X_i\|_{r_i})+\\
&+(\E (\sum_{i\in I(t)\bigtriangleup I(s)} k_i |X_i|1_{|X_i|\gs \alpha e\|X_i\|_{r_i}})^p )^{\frac{1}{p}}\ls \\
&\ls  (\alpha e \sum_{i\in I(t)\bigtriangleup I(s)} k_i\|X_i\|_{r_i})+ \alpha \gamma e(\E |\sum_{i\in I(t)} |\ccE_i|1_{|\ccE_i|\gs r_i}|^p )^{\frac{1}{p}}\ls \\
&\ls (\alpha e \sum_{i\in I(t)\bigtriangleup I(s)} k_i \|X_i\|_{r_i})+\alpha \gamma e \|Z\|_p,
\end{align*}
where $Z$ is of gamma distribution $\Gamma(|I(t)|,1)$. 
Clearly $\|Z\|_p\ls (p+|I(t)|)\ls (1+c)p\ls 2p$, and hence for $\gamma= (8\alpha e)^{-1}$
$$
\alpha \gamma e \|Z\|_p \ls \frac{p}{4}.
$$
Consequently using that $\|X_t-X_s\|_p\gs \frac{p}{2}$ it implies that
$$
\frac{p}{4}\ls \alpha e \sum_{i\in I(t)\bigtriangleup I(s)} k_i \|X_i\|_{r_i}\ls \alpha e \sum_{i\in I(t)\bigtriangleup I(s)}r_i.
$$
Finally for sufficiently small $c$ we have $2\alpha e|I(t)\bigtriangleup I(s)|\ls 4c\alpha e p\ls \frac{p}{8}$, so
$$
\frac{p}{8}\ls \alpha e\sum_{i\in I(t)\bigtriangleup I(s)}r_i 1_{r_i> 2}
$$
which ends the proof of the result. 
\end{dwd}
We turn to prove the Sudakov minoration for independent $X_1,...,X_n$. 
\smallskip

\noindent
The first tool we apply is the reduction of the problem to just 
symmetric independent exponentials. This is due to the Bernoulli comparison
in its most powerful form. Note that here we need the lower tail bound (\ref{r2}) for each $X_i$,
$i\in \{1,2,...,n\}$.
\begin{lema}\label{lema2}
Let $(Y_1,...,Y_n)$ be independent symmetric and sub-exponential in the following  sense
$$
\P(|Y_i|>u)=e^{-u},\;\;u\in [0,r_i].
$$
Then
$$
\E \sup_{t\in T}|X_t|\gs \frac{1}{C^2}\E \sup_{t\in T}|\sum_{i\in I(t)} \frac{\|X_i\|_{r_i}}{r_i} k_i Y_i|.
$$
\end{lema}
\begin{dwd}
The so called Bernoulli comparison (see Lemma 4.6 in \cite{Le-Ta}) states that
for two sequences of independent symmetric variables $\eta_i,\xi_i$, $i\in \{1,2,...,n\}$ 
the comparability of tails
\be\label{mandy}
C\P(|\xi_i|>u)\gs \P(|\eta_i|>u),\;\;\mbox{for all}\;u\gs 0,\;\;\mbox{and}\;\;i\in \{1,2,...,n\}
\ee
for two independent, symmetric families of variables $\eta_i, \xi_i$, $i\in \{1,2,...,n\}$ then
\be\label{andy}
C\E\sup_{t\in T}|\sum^n_{i=1} t_i \xi_i| \gs \E\sup_{t\in T}|\sum^n_{i=1} t_i \eta_i|.
\ee
We apply the result for $\eta_i=\ccE_i$ and $\xi_i=\frac{r_i X_i}{C\|X_i\|_{r_i}}$. Obviously (\ref{r2})
implies (\ref{mandy}) and hence (\ref{andy}) holds which is the acquired inequality. 
\end{dwd}
The second tool is  the basic minoration
for exponentials \cite{Lat}. Note that this is result is the core of the proof
and is based on a tricky induction which is difficult to repeat for the dependent case.
\begin{prop}\label{prop4}
Suppose that for any $s,t\in T$
$$
\sum_{i\in I(t)\backslash I(s)} r_i v_i \gs q,
$$
where $r_i,v_i\gs 1$ and $T\gs e^q$ then
$$
\E \sup_{t\in T} |\sum^n_{i=1}v_i Y_i|\gs \frac{1}{8}q.
$$
\end{prop}
We have collected all the tools to complete the proof of the Sudakov minoration for independent entries. The result was announced in \cite{La-Tk} we give our proof for the sake of completeness
of this paper.
\begin{theo}\label{theo1}
There exists a universal constant $K$ such that for any set $T$ of the form stated in Proposition \ref{prop1}
the following inequality holds
$$
\E\sup_{t\in T}|X_t|\gs K^{-1}p.
$$
\end{theo}
\begin{dwd}
By the Proposition \ref{prop1} we have 
$$
\| X_t-X_s\|_{p}\gs \frac{p}{2},\;\;\mbox{for all}\;s,t\in T,\;s\neq t.
$$
Lemma \ref{lema3} implies that for all $s,t\in T$, $s\neq t$
\be\label{bum}
\sum_{i\in I(t)\bigtriangleup I(s)}r_i1_{r_i> 2}\gs \gamma p.
\ee
Then by Lemma \ref{lema2} and the definition of $r_i$ we obtain that
\begin{align*}
& \E \sup_{t\in T}|X_t|\gs C^{-2}\E \sup_{t\in T}|\sum_{i\in I(t)} \frac{\|X_i\|_{r_i}}{r_i}k_iY_i|\gs\\
&\gs  C^{-2}\E |\sup_{t\in T}\sum_{i\in I(t)}1_{r_i > 2}Y_i|.
\end{align*}
We aim to apply Proposition \ref{prop4}. Let $q=\min\{4^{-1},\gamma\}p$.
Let $v_i=1_{r_i> 2}$. Then by (\ref{bum})
$$
\sum_{i\in I(t)\bigtriangleup I(s)} r_i v_i\gs q.
$$
Consequently by Proposition \ref{prop4}
$$
\E \sup_{t\in T} |\sum_{i\in I} 1_{r_i>2} Y_i|\gs \frac{q}{8}
$$  
and hence for a universal $K$
$$
\E \sup_{t\in T}|X_t|\gs K^{-1}p,
$$
which completes the proof.
\end{dwd}

\section{Disjoint supports}\label{sect4}

The next step is to generalize the Sudakov type minoration on the cases where there is some dependence among the entries of $X$.
From now on we assume that $X=(X_1,...,X_n)$ is one unconditional and log concave.
The simplest case in which one can try to prove Sudakov minoration concerns $T$ where all the points are of disjoint support
i.e. $I(t)\cap I(s)=\emptyset$ for all $s,t\in T$, $s\neq t$. We give our proof that in this setting the Sudakov minoration indeed
works and then deduce from it that function $f(p)=p^2$ is right upper bound on the cardinality of $T$ that implies
Sudakov minoration. The proof is based on (\ref{n6}) - the exponential inequality for log-concave distribution. 
Note that the result does not hold without regularity assumptions on the one unconditional distribution.
\begin{theo}\label{theo2}
Suppose that $|T|\gs \exp(C p)$, for $C\gs 1$, sufficiently large. Suppose that $0\in T$ and all points $t\in T$ have disjoint supports and 
\be\label{n7}
\| X_t\|_p\gs p,\;\;\mbox{for all}\;t\in T,\;t\neq 0
\ee
then there exists a universal constant $K$ such that
$$
\E \sup_{t\in T}|X_t|\gs K^{-1}p.
$$
\end{theo}
\begin{dwd}
First observe that $\langle t, X\rangle /\|t\|_2$, $t\in T$, $t\neq 0$ is still isotropic and log concave vector. 
Enumerate points in $T$ as $t_0,t_1,...,t_N$, where $t_0=0$ and obviously $N=|T|-1$. Denote $a_i=\|t_i\|_2$, $Y_i=\langle t_i,X\rangle/\|t_i\|_2$ for
$i\in \{1,...,N\}$. The assumption (\ref{n7}) implies that
\be\label{mradek1}
a_i \|Y_i\|_p\gs p,\;\;\mbox{for all}\; 1\ls i\ls N.
\ee
For all $i\in \{1,2,...,N\}$ let
$$
S_i=\{x\in \R^N: a_i|x_i|\gs \bar{C}^{-1}_0 p\},
$$
where $\bar{C}_0\gs 1$ is a universal constant such that
$$
\P(Y\in S_i)= \P(a_i|Y_i|\gs \bar{C}^{-1}_0 p)\gs 2e^{-p}.
$$
The existenceof such $\bar{C}_0$ can be deduced from (\ref{mradek1}), (\ref{r2}) and the convexity of $u\ra -\log \P(|Y_i|\gs u)$.
We count how many variables $a_i|Y_i|$, $i\in \{1,2,...,N\}$ crosses the level $\bar{C}^{-1}_0 p$, i.e.
we introduce the following variable
$$
M=\sum^N_{i=1}1_{a_i|Y_i|\gs \bar{C}^{-1}_0 p}. 
$$
Clearly
$$
N^{-1}\E M =N^{-1}\sum^N_{i=1} \P(a_i|Y_i|\gs \bar{C}^{-1}_0 p)\gs  2e^{-p}.
$$
Now we choose $\bar{C}_1\gs 2\bar{C}_0$. We can assume that 
\be\label{n10}
\P(\max_{1\ls i\ls N}a_i|Y_i|\ls \bar{C}_1^{-1}p )\gs \frac{1}{2}.
\ee
Indeed suppose that $\P(\max_{1\ls i\ls N}a_i|Y_i|>\bar{C}_1^{-1}p)\gs \frac{1}{2}$, then
$$
\E \max_{1\ls i\ls N} a_i |Y_i|\gs \frac{1}{2}\bar{C}^{-1}_1 p.
$$
which is the acquired minoration.
Observe that we can always assume $a_i\gs \bar{C}_2^{-1}p$ for a universal $\bar{C}_2$.
Otherwise
$$
\E \max_{1\ls i\ls N} a_i |Y_i|\gs \bar{C}_2^{-1}p\E |Y_j|
$$
for some $j\in \{1,2,...,N\}$ and $\E |Y_j|\gs \bar{C}_3$ due to the isotropy nad log concavity of $Y$. 
Thus again there would be nothing to prove. Let
$$
S=\{x\in \R^N:\; a_i|x_i|\ls \bar{C}_1^{-1}p\}.
$$
Clearly (\ref{n10}) means $\P(X\in S)\gs \frac{1}{2}$.
\smallskip

\noindent
We choose $n_0= e^{-p}N$ which guarantees that
$$
2e^{-p}\ls N^{-1}\E M\ls e^{-p}+N^{-1}\E M 1_{M>n_0}\ls e^{-p}+\P(M>n_0),
$$
and therefore
\be\label{n9}
\P(M>n_0)\gs e^{-p}.
\ee
We have to understand the geometry of the set $\{M\gs k\}$ for $k> n_0$.
It means that there exists set $K$ of large cardinality such that $a_i|Y_i|\gs \bar{C}^{-1}_0 p$
for all $i\in K$, namely
$$
\{M\gs k\}=\{\exists K\subset \{1,2,...,N\}:\;|K|=k\;a_i|Y_i|> \bar{C}^{-1}_0 p,\;\forall\; i\in K\}.
$$
Fix $K\subset \{1,2,...,N\}$, $|K|=k$. Consider set $S_K$ of the form 
$$
S_K=\bigcap_{i\in K}S_i=\{x\in \R^n:\; a_i|x_i|> \bar{C}^{-1}_0 p,\; \forall \; i\in K \}.
$$
We show that $S_K$ is well separated from $S$ in the sense of $\ell^2$ distance, i.e.
\be\label{n8}
S+\beta u B^2_N\cap S_K=\emptyset,\;\;\mbox{where}\;u=\bar{C}_3^{-1}\beta^{-1}k^{\frac{1}{2}}
\ee
for a universal constant $\bar{C}_3$. Consider point $z=x-y$, where $x\in S_K$ and $y\in S$.
If $x\in S_K$ then $a_i|x_i|\gs \bar{C}^{-1}_0 p$ for all $i\in K$. 
On the other hand if $y\in S$ then $a_i|y_i|\ls \bar{C}_1^{-1}p$,  for all $i\in \{1,2,...,N\}$.  Hence
$$
a_i|z_i|\gs a_i(|x_i|-|y_i|)\gs (\bar{C}^{-1}_0-\bar{C}^{-1}_1)p\gs \bar{C}_0^{-1}\frac{p}{2},\;\;\mbox{for all}\;i\in K,
$$ 
and therefore
$$
\|z\|^2_2\gs 2^{-2}(\bar{C}_0)^{-2} p^2 \sum_{i\in K} (a_i)^{-2}.
$$
As we have mentioned $a_i\ls \bar{C}_2^{-1} p$ and consequently for any 
$K\subset \{1,2,...,N\}$ such that  $|K|=k$
$$
\|z\|_2> \bar{C}_3^{-1} k^{\frac{1}{2}},
$$
where $\bar{C}_3=2\bar{C}_0\bar{C}_2$. It proves (\ref{n8}). 
\smallskip

\noindent
Since it works for any $K$, $|K|\gs k>n_0$
we can apply our main tool i.e. the exponential inequality  (\ref{n6}) which gives 
$$
\P(M>n_0)=\P(\exists K, |K|>n_0:\;Y\in S_K)\ls e^{-u},\;\;\mbox{for}\;u=\bar{C}_3 \beta^{-1}n_0^{\frac{1}{2}}.
$$
Obviously we need that $u>p$ which means
\be\label{anatol}
\bar{C}_3 n_0^{\frac{1}{2}}=\bar{C}_3(e^{-p}N)^{\frac{1}{2}}> \beta.
\ee
This is the point where our main assumption on $\beta$ matters. Indeed $\beta\ls \bar{C}_4 N^{\frac{1}{2}-\va}$ implies that
for $N=\exp(C p)-1$ with $C$ large enough we can compensate the value of $\va$ and guarantee that (\ref{anatol}) holds. 
In this way we end up with contradiction
$$
e^{-p}\ls N^{-1}\E M\ls \frac{n_0}{N}+P(M>n_0)< e^{-p}. 
$$ 
Consequently what fails is the assumption that $\P(Y\in S)\gs \frac{1}{2}$ and hence the proof is completed.
\end{dwd}
The basic consequence of the above result is that the Sudakov minoration works for log concave one unconditional $X$
whenever $|T|\gs \exp(Cp^2)$, (i.e. $f(p)=Cp^2$) and $C$ is arbitrary small, with the constant $K$ that depends on $C$ only. We improve slightly the result using the
combinatorial dimension of the class $I(t)$, $t\in T$. Recall that class $\ccC$ of subsets of $\{1,2,...,n\}$ is of dimension $v$ if
there is no subset of $\{1,2,...,n\}$ of cardinality $v+1$ that is shattered by the class $\ccC$ into all its subsets. Since
$|I(t)|\ls cp$ where $c$ can be sufficiently small we have that at least $v\ls cp$.
\begin{coro}\label{cor4}
Suppose that $I(t)$, $t\in T$ has a VC dimension $v\ls cp$ and $|T|\gs \exp(f(p))$ for $f(p)=(c+C)vp$, where the constant $C$ is from Theorem \ref{theo2} then
whenever $T$ is of the form stated in Proposition \ref{prop2} with $q=p$ then following inequality holds
$$
\E \sup_{t\in T} X_t\gs K^{-1}p,
$$
where $K$ is a universal constant.
\end{coro}
\begin{dwd}
We simply show that after the simplification from Proposition \ref{prop2}
and $f(p)=v(C p+1)$ there must exist at least $\exp(C p)$ points of disjoint supports. 
\smallskip

\noindent
First apply Proposition \ref{prop2}. Recall that by Corollary \ref{cor2} all the supports $I(t)$, $t\in T$ are thin $|I(t)|\ls cp$
and different. We start the construction from $t^0=0$ and $I_0=J_0=\emptyset$ and then continue construction of $I_l$, $J_l$ and $t_l$ in the following way.
Suppose that the set 
$$
S_{l+1}=\{t\in T:\; \|X_{t}-X_{t}1_{J_l}\|\gs \frac{p}{8}\}
$$ 
is not empty. We select $t^{l+1}$ as any point in $S_{l+1}$ and $I_{l+1}=I(t^{l+1})\backslash J_l$
and $J_{l+1}=J_l\cup I(t^{l+1})$. The construction stops after $N$ steps, which means $S_{N+1}$ is empty whereas $S_N$ is not.
Clearly points $s^l=t^l 1_{I_l}$ are of disjoint supports and 
$$
\|X_{s^l}\|_p\gs \frac{p}{8}\;\;\mbox{for all}\; 1\ls l\ls M.
$$
Therefore if $N\gs \exp(Cp)$ we deduce 
$$
\E \sup_{0\ls l\ls N}X_{s^l}\gs K^{-1}p
$$ 
and hence due to the Bernoulli comparison
$$
\E \sup_{t\in T} X_t\gs\E \sup_{0\ls l\ls M}X_{s^l}\gs  K^{-1}p.
$$
Now if $M<\exp(Cp)$ then we can consider points $t 1_{J_N}$ and observe that
$$
\|X_{t}1_{J_N} -X_s 1_{J_N}\|_p\gs \frac{p}{4}\;\;\mbox{for all}\;s,t\in T,s\neq t.
$$
Indeed $\|X_t-X_s\|_p\gs \frac{p}{2}$ implies that
\begin{align*}
& \|X_{t}1_{J_N} -X_s 1_{J_N}\|_p\gs \|X_t-X_s\|_p- \|X_t-X_{t}1_{J_N}\|_p-\|X_s-X_{s}1_{J_N}\|_p\gs \\
&\gs \frac{p}{2}-\frac{p}{8}-\frac{p}{8}=\frac{p}{4}.
\end{align*}
It means that points $t 1_{J_N}$ for $t\in T$ are well separated and due to the same argument as in Corollary \ref{cor2}
also of different supports. The set $J_N$ counts no more then $cpN$ points which is small than $\exp((c+C)p)$. We use Sauer's lemma (e.g Proposition 14.10 in \cite{Le-Ta}).
Let $\bar{C}=1+c+C$ there is less than
$$
\left(\frac{e^{\bar{C}p}}{v}\right)^{v}\ls \exp(\bar{C}vp)
 $$
elements possible in $J_M$ of supports of cardinality not larger than $v$. Therefore we have a contradiction if $|T|\gs \exp(\bar{C}vp)$ which completes the proof.
\end{dwd}

\section{Common witness}\label{sect5}

Now we turn to prove some extension of Theorem \ref{theo2}. Our aim is to slightly relax the assumption that
supports are disjoint and prove that $|T|\gs \exp(C p\log(1+p))$ suffices if for each point $t\in T$
there exists a common witness. Note that $f(p)=Cp\log(1+p)$ is much better than $p^2$ 
we have proved to be universal bound for log-concave unconditional vectors.
\smallskip

\noindent 
Recall that by Proposition \ref{prop2} and consequently by Corollary \ref{cor2} we have that supports of $t\in T$
are thin ($|I(t)|\ls cp$) and different ($I(t)\neq I(s)$ for $s,t\in T$, $s\neq t$).
By Corollary \ref{cor3} the condition $\|X_t-X_s\|_p\gs \frac{p}{2}$ implies that 
\be\label{trach}
\|\sum_{i\in I(t) \backslash I(s)}k_i 1_{k_i> 4\rho} X_i\|_p\gs \frac{p}{8}\;\;\mbox{or}\;\;\|\sum_{i\in I(s) \backslash I(t)}k_i 1_{k_i> 4\rho} X_i\|_p\gs \frac{p}{8}.
\ee 
We need that $k_i> 4\rho$ since it guarantees that $t_i\ls k_i\ls 2t_i$.
Fix $t\in T$. Define set $S(t)\subset T$ of significant neighbors of $t$ by 
$$
S(t)=\{s\in T: \|\sum_{i\in I(t)\backslash I(s)} k_i 1_{k_i> 4\rho} X_i\|_p \gs \frac{p}{8}\}.
$$
By (\ref{trach}) we get that either $s\in S(t)$ or $t\in S(s)$. 
Theorem \ref{bambam} implies that whenever $s\in S(t)$ there exists point $a(t,s)\in R^n$ such that 
$$
\P(\bigcap_{i\in I(t)\backslash I(s)}\{|X_i|\gs a_i(t,s)\})\gs e^{-p}
$$
and $\sum_{i\in I(t)\backslash I(s)}k_i 1_{k_i> 4\rho} a_i(t,s)\gs D^{-1}\frac{p}{8}$. 
Our basic assumption in this section is that for all $t\in T$ there exist a common witness $a(t)$, i.e.
$a(t)\in \R^n$ such that 
\be\label{panama1}
\P(\bigcap_{i\in I(t)}\{|X_i|\gs a_i(t)\})\gs e^{-p}
\ee
and for all $s\in S(t)$ there holds
\be\label{panama2}
\sum_{i\in I(t)\backslash I(s)} k_i 1_{k_i> 4\rho} a_i(t)\gs C^{-1}_4p
\ee
for a universal constant $C_4$. Note that we may easily require that $a_i(t)\gs C^{-1}_5$, where $C_5$
is an absolute constant. One can either deduce it straight from Theorem \ref{bambam1} or use the following argument.
Due to (\ref{brum1}) we have
$$
C^{-1}_5\sum_{i\in I(t)} k_i\ls 2C_0C^{-1}_5\delta p\ls C^{-1}_4 \frac{p}{2}
$$
for small enough $\delta$. Hence we can still have the lower bound (\ref{panama2}) when $a_i(t)\gs C^{-1}_5$.
By Theorem \ref{bambam1} the existence of a common witness for $t\in T$ is equivalent to
\be\label{panama3}
\|\min_{s\in S(t)}|\sum_{i\in I(t)\backslash I(s)} k_i1_{k_i> 4\rho} X_i|\|_p\gs C^{-1}_6 p.
\ee
There are many cases where the condition holds we list some of them.
\begin{enumerate}
\item  Disjoint supports. Obviously $S(t)=T\backslash \{t\}$ and $I(t)\backslash I(s)=I(t)$, so the existence of a 
common witness is the same as the existence of a witness for $t\in T$.
\item  Domination of supports. In this case we assume that
$\|\sum_{i\in I(t)} k_i 1_{k_i> 4\rho}X_i\|_p$ slightly dominates overlaps i.e. there exists $\va>0$ such that for all $s\in S(t)$
\be\label{puzian}
\|\sum_{i\in I(t)\cap I(s)}k_i 1_{k_i> 4\rho} X_i\|_p\ls (1-\va)\|\sum_{i\in I(t)} k_i 1_{k_i> 4\rho} X_i\|_p.
\ee
We prove that if we choose $c$ to be sufficiently small such that $(1-\va)2^c<1$ then (\ref{puzian}) implies the existence
of a common witness.
\begin{lema}
Suppose that there exits $\va>0$ such that $(1-\va)2^u<1$ and 
$$
\|\sum_{i\in I(t)\cap I(s)}k_i1_{k_i> 4\rho} X_i\|_p\ls (1-\va)\|\sum_{i\in I(t)} k_i1_{k_i> 4\rho} X_i\|_p\|\;\;\mbox{for all}\;\;s\in S(t).
$$
Then
$$
\|\min_{s\in S(t)}|\sum_{i\in I(t)\cap I(s)}k_i1_{k_i> 4\rho} X_i\|_p\gs (1-(1-\va)2^c)p.
$$
\end{lema}
\begin{dwd}
We aim to prove that (\ref{panama3}). Clearly
\begin{align*}
& \|\min_{s\in S(t)}|\sum_{i\in I(t)\backslash I(s)} k_i1_{k_i> 4\rho} X_i|\|_p\gs\\
&\gs \|\sum_{i\in I(t)} k_i1_{k_i> 4\rho} X_i\|_p-\|\max_{s\in S(t)}\sum_{i\in I(t)\cap I(s)}k_i1_{k_i> 4\rho} X_i\|_p\gs \\
&\gs \|\sum_{i\in I(t)} k_i1_{k_i> 4\rho} X_i\|_p-(\sum_{s\in S(t)}\E |\sum_{i\in I(t)\cap I(s)}k_i1_{k_i> 4\rho} X_i|^p )^{\frac{1}{p}}\gs\\
&\gs (1-(1-\va)2^c)\|\sum_{i\in I(t)} k_i1_{k_i> 4\rho} X_i\|_p.
\end{align*}
\end{dwd}
\item Independent entries. This case is of particular interest since there is a need for a different proof
of the Sudakov minoration than the induction argument we have presented in section \ref{sect3}. Recall the definition of $r_i$, $1\ls i\ls n$ and Lemma \ref{lema3}.
It is clear that there exists the following distance on $T$
$$
\bar{d}(s,t)=\sum_{i\in I(t)\bigtriangleup I(s)} r_i 1_{r_i> 2}. 
$$ 
The best setting for our purposes is when $\bar{d}(s,t)$ is of finite distortion, i.e. when 
$$
\bar{C}^{-1}q\ls \bar{d}(s,t)\ls \bar{C}q\;\;\mbox{for}\;q\gs p.
$$
This requires a slight generalization of the common witness definition. We say that $q$ common witness exists
for $t\in T$ and $q\gs p$ if there is $a(t)\in \R^n$ supported in $I(t)$ such that   
$$
\P(\bigcap_{i\in I(t)}\{|X_i|\gs a_i(t)\})=e^{-q}
$$
and for all $s\in S_q(t)$, where 
$$
S_q(t)=\{s\in T:\; \sum_{i\in I(t)\backslash I(s)}r_i 1_{r_i> 2}\gs 2^{-1}\bar{C}^{-1}q\}
$$
the following holds
$$
\sum_{i\in I(t)\backslash I(s)} k_i1_{k_i> 4\rho} a_i(t)\gs \bar{C}_0^{-1} q.
$$ 
First observe that for each $s,t\in T$ the condition $\bar{d}_q(s,t)\gs \bar{C}q$  implies that either $s\in S_q(t)$ or $t\in S_q(s)$. 
By the definition of $r_i$ it implies that for all $s\in S_q(t)$
$$
\bar{C}^{-1}_0\frac{q}{2}\ls \sum_{i\in I(t)\backslash I(s)} r_i 1_{r_i> 2}\ls \gamma^{-1}\sum_{i\in  I(t)\backslash I(s)} k_i\|X_i\|_{r_i} 1_{r_i> 2}.
$$
Consequently for $a_i(t)=\bar{C}^{-1}_1\|X_i\|_{r_i}1_{r_i> 2}$, $i\in I(t)$ where $\bar{C}_1$ is an absolute  constant we have
$$
\sum_{i\in I(t)\backslash I(s)} k_ia_i(t)\gs \bar{C}^{-1}_0\bar{C}_1^{-1}\gamma \frac{q}{2}.
$$
Suppose that $k_i\ls 4\rho$ then $k_i\|X_i\|_2=k_i\ls 4\rho\ls 2\gamma$ for sufficiently small $\rho$. Hence
$k_i\ls 4\rho$ implies that $r_i=2$ and therefore
$$
\sum_{i\in I(t)\backslash I(s)} k_ia_i(t)=\sum_{i\in I(t)\backslash I(s)} k_i1_{k_i> 4\rho}a_i(t).
$$
On the other hand using log-concavity and suitably choosing $\bar{C}_1$
\begin{align*}
& \P(\bigcap_{i\in I(t)}\{|X_i|\gs a_i(t)\})=\prod_{i\in I(t)}\P(|X_i|\gs \bar{C}^{-1}_1 \|X_i\|_{r_i}1_{r_i> 2})\gs \\
& \gs \prod_{i\in I(t)}\exp(-\bar{C}^{-1}_0 r_i1_{r_i>2})=
\exp(-\bar{C}_0^{-1}\sum_{i\in I(t)}r_i1_{r_i>2})\gs e^{-q}.
\end{align*}
In this way having the distortion controlled we are in the setting of a common witness existence.
Since $\bar{d}$ is a distance we can always find a suitably large set of 
bounded distortion for some large enough $q$ if we lose slightly on the power of the function $f$.
Namely the following holds.
\begin{lema}
Suppose that $|T|\gs \exp(\bar{f}(p))$ where $\bar{f}(p)=\log(1+p)f(p)$ then there exists $S\subset T$
such that $|S|\gs \exp(f(p))$ such that $C^{-1}q\ls \bar{d}(s,t)\ls Cq$ for all $s\neq t$, $s,t\in S$, where $C$
is a universal constant.
\end{lema}
\begin{dwd}
We use the chaining type argument to get this result. The crucial is to understand that $\bar{d}(s,t)\ls 2cp^2$
since $r_i\ls p$ and $|I(s)|,|I(t)|\ls cp$. Therefore $M=\sup_{s,t}\bar{d}(s,t)\ls cp^2$. On the other hand 
by Lemma \ref{lema3} we have that $\bar{d}(s,t)\gs \bar{C}_0^{-1}p$ for some constant $p$.
\smallskip

\noindent
Suppose that there exists $S\subset T$ such that $|S|\gs \exp(f(p))$ and $\bar{d}(s,t)\gs \bar{C}_1^{-1}M$
for all $s,t\in S$ and $s\neq t$. Then simply we put $q=M$ and finish the construction.
Otherwise there exists a partition of $T$ into $T_1,...,T_N$ where $N\ls \exp(f(p))$ such that
$\sup_{s,t\in T_i}\bar{d}(s,t)\ls 2\bar{C}_1^{-1}M$. We continue the construction in this way.
Consequently since $M\ls cp^2$ after $\bar{C}_2^{-1}\log(1+p)$ steps, where $\bar{C}_2$ can be sufficiently large
we either find set $S$ and $p\ls q\ls M $ such that $\bar{C}^{-1}_1 q\ls \bar{d}(s,t)\ls q$ for all $s,t\in S$
or we end up with less or equal $\exp(\bar{C}_2^{-1}\log(1+ p)f(p))$ sets $T_1,...,T_N$ that covers $T$
such that $\bar{d}(s,t)\ls p$. However since $|T|\gs \exp(\log(1+p) f(p))$ it means that at least one set
$T_i$ counts at least $\exp(f(p))$ points since otherwise 
$$
\sum^N_{i=1}|T_i|\ls \exp(\bar{C}_2^{-1}\log(1+ p)f(p))\exp(f(p))< \exp(\log(1+p)f(p))
$$
for suitably large. This contradiction completes the proof. 
\end{dwd}

\end{enumerate}
Before we state the main result we explain the meaning of a common witness. 
We will be able to prove that whenever a common witness exists the following holds
$$
\E \sup_{t\in T} \sup_{s\in S(t)}|X_t 1_{I(t)\backslash I(s)}|\gs K^{-1}p,
$$
where $K$ is an absolute constant. As it was proved in Lemma \ref{lema5} in the ideal case
it is equivalent to the Sudakov minoration. Sometimes it is obvious by the
conditions imposed on points in $t\in T$ that guarantees   
$$
C \E \sup_{t\in T} |X_t|\gs \E \sup_{t\in T} \sup_{s\in S(t)}|X_t 1_{I(t)\backslash I(s)}|
$$
for a universal constant $C$. This for example the case of disjoint supports where $I(s)$ is empty for all $s\in S(t)$.
In general it is only true that
$$
\E \sup_{t \in T}|X|_t \gs \E \sup_{t\in T} \sup_{s\in S(t)}|X_t 1_{I(t)\backslash I(s)}|.
$$
The Sudakov minoration $\E \sup_{t\in T}|X|_t\gs K^{-1}p$ is considered to be of the similar difficult
as the standard Sudakov minoration. Note that as we have mentioned in Corollary \ref{cor3}
after the simplification from Proposition \ref{prop2} the one unconditional structure of $X$ does not matter for the
assumption on increments.
\smallskip

\noindent
We are in the position to prove an extension of Theorem \ref{theo2} which is the main new result of the paper.
\begin{theo}
Suppose that after simplification from Proposition \ref{prop2} the class of $I(t)$, $t\in T$ is of VC dimension 
$v$ ($v\ls cp$). Suppose that $|T|\gs \exp(f(p))$, where $f(p)=C v\log(1+p)$ where $C$ is sufficiently large. 
Suppose that for each $t\in T$ the common witness exists, i.e.
$$
\|\min_{s\in S(t)}|\sum_{i\in I(t)\backslash I(s)}k_i1_{k_i>4\rho} X_i|\|_p\gs \bar{C}^{-1}p.
$$
Then the following Sudakov minoration holds
$$
\E \sup_{t\in T}\sup_{s\in S(t)}|X_{t}1_{I(t)\backslash I(s)}|\gs K^{-1}p
$$
for a universal $K$. In particular if $v\ls \frac{p}{\log(1+p)}$  then $f(p)=Cp$.
\end{theo}
\begin{dwd}
First we assume that $n\ls \exp(\bar{C}_0 p)$ since otherwise we may apply Theorem \ref{theo2} and there is nothing to 
prove. We need this condition to control the constant $\beta$ in the exponential inequality by $\bar{C}_1p$.  
\smallskip

\noindent
Let us enumerate points in $T$ by $t^0,t^1,...,t^N$, where $t^0=0$ and $N=|T|-1$. 
As in Theorem \ref{theo2} we aim to contradict the assumption that
$$
\E\sup_{0\ls i\ls N}|X_{t^i}1_{I(t^i)\backslash I(t^j)}|\ls \bar{C}^{-1}_2 \frac{p}{2} 
$$
for some absolute $\bar{C}_2$. Consequently 
\be\label{koyot1}
\P(\sup_{1\ls i\ls N}|X_{t^i} 1_{I(t^i)\backslash I(t^j)}|\ls \bar{C}^{-1}_2 p)\gs \frac{1}{2}.
\ee
Let $S(i)=S(t_i)$ and
$$
S=\{x\in \R^n:\; |\langle t^i 1_{I(t^i)\backslash I(t^j)}, x\rangle|  \ls \bar{C}_2^{-1} p,\;\;\mbox{for all}\;i\in \{1,2,...,N\}, j\in S(i)\}.
$$
Note that for all $x\in S$ for all $j\in S(i)$ and $i\in \{1,2,...,N\}$
\be\label{onion}
\sum_{l\in I(t_i)\backslash I(t_j)} t^i_l x_l\ls \bar{C}^{-1}_2p. 
\ee
The meaning of (\ref{koyot1}) is that $\P(X\in S)\gs \frac{1}{2}$ and the set $S$ will play the
same role in the proof as in Theorem \ref{theo2}.
With each point $i\in \{1,2,...,N\}$ we choose a common witness $a^i=a(t_i)$ and hence select the
set
$$
S(i)=\{x\in\R^n:\; x_j\gs \bar{C}^{-1}_3a^i_j,\; j\in I(t_i)\},
$$
where $\bar{C}_3\gs 1$. By the definition of a common witness
$$
\P(X\in S(i))\gs 2^{-cp}\P(\bigcap_{j\in I(t_i)} \{|X_j|\gs a^i_j\} )^{\bar{C}_3^{-1}}\gs  \exp(-(c+\bar{C}_3^{-1})p).
$$
Choosing $c$ sufficiently small and $\bar{C}_3$ large enough (say $\bar{C}_3=4$) we can guarantee that 
$$
\P(X\in S(i))\gs 2e^{-p}.
$$
As in Theorem \ref{theo2} let us define  $M=\sum^N_{i=1}1_{X\in S(i)}$, so
\be\label{bekon}
\frac{1}{N}\E M\gs 2e^{-p}. 
\ee
Since $N\gs \exp(f(p))-1=\exp(\bar{C}_{0}v\log(1+p))-1$ where we may require $v$ to be greater or equal $\frac{p}{\log (1+p)}$
we deduce that defining $n_0=e^{-p}N$ we have that
$$
\frac{1}{N}\E M\ls \frac{n_0}{N}+\P(M\gs n_0)\ls e^{-p}+\P(M\gs n_0)
$$
and hence $\P(M\gs n_0)\gs e^{-p}$ by (\ref{bekon}). 
\smallskip

\noindent
In the view of the idea of the proof of Theorem \ref{theo2} we have to consider sets
$S_K=\bigcap_{i\in K}S_i$, where $|K|\gs e^{-p}N$ . We aim to show that $S+uB^n_2$ does not intersect $S_{K}$ for a sufficiently large $c$, namely
that there exists sufficiently large constant $\bar{C}_{4}$ such that
\be\label{knot}
S+\beta uB^n_2\cap S_K=\emptyset,\;\;\mbox{for}\; u\ls \beta^{-1}p^{\bar{C}_{4}}.
\ee
Consider point $z=x-y$, where $x\in S_K$ and $y\in S$. We prove that there exists at least $p^{\bar{C}_{5}}$
coordinates $l\in \{1,2,...,n\}$ such that $|z_l|\gs \bar{C}_{6}$ and $\bar{C}_{5}$ and $\bar{C}_{6}$ are sufficiently large.
\smallskip

\noindent
We start from $L_1=\{i\}$ a single $i\in K$. There must exists a single $l(i)\in I(t)$ such that $y_{l(i)}\ls \bar{C}^{-1}_7 a^i_{l(i)}$, where $\bar{C}_7$
is sufficiently small. Otherwise $y_l\gs C^{-1}_7 a^i_{l}$ for all $l\in I(t^i)$ and hence by (\ref{brum1}), (\ref{panama2}) and $t^i_l\ls k_l\ls 2t^i_l$ for $k_l\gs 4\rho$
$$
\sum_{l\in I(t^i)} t^i_l y_l\gs 2^{-1}\bar{C}^{-1}_7\sum_{l\in I(t^i)} k_l 1_{k_l> 4\rho} a^i_l \gs 2^{-1}C_4^{-1}\bar{C}^{-1}_7 p.  
$$
Therefore if $\bar{C}_2> 2C_4 \bar{C}_7$ we have a contradiction with (\ref{onion}).  Suppose we have selected set $L_k\subset K$ such that $|L_k|=k$
and for each $i\in L_k$ there exists at least one $l(i)\in I(t^i)$ such that $y_{l(i)}\ls \bar{C}_7^{-1}a^i_{l(i)}$. We require that 
coordinates $l(i)$, $i\in L_k$ are all different. Consider the set $J=\bigcup_{i\in L_k} I(t^i)$. By our basic inequality
$|I(t^i)|\ls cp$ we deduce that $|J|$ counts no more than $cpk$ elements. 
Therefore by the Sauer's lemma (e.g. \cite{Le-Ta} chapter 14) we can create at most $(\frac{ckpe}{v})^v$ possible 
subsets of $J$. Consequently if $K$ has more elements that $\left(\frac{ckpe}{v}\right)^v$ there are two points $t^i,t^j$, $i\neq j$, $i,j\in K$
such that $I(t_i)\cap J=I(t_j)\cap J$ and hence $(I(t_i)\bigtriangleup I(t_j))\cap J=\emptyset$.
Therefore either for $i$ or $j$ there must exist $l(i)$ or $l(j)$ such that
$$
y_{l(i)}\ls \bar{C}_7^{-1}a^i_{l(i)}\;\;\mbox{or}\;\;y_{l(j)}\ls 2^{-1}\bar{C}_7^{-1}a^j_{l(j)}
$$ 
since otherwise $y_l\ls \bar{C}_7^{-1}a^i_l$ for all $l\in I(t^i)\backslash I(t^j)$
and $y_l\ls \bar{C}_7^{-1}a^j_l$ for all $l\in l\in I(t^j)\backslash I(t^i)$ which implies
$$
\sum_{l\in I(t^i)\backslash I(t^j)} t^i_l1_{k_l> 4\rho} y_l\gs 2^{-1}\bar{C}^{-1}_7\sum_{l\in I(t^i)\backslash I(t^j)} k_l1_{k_l> 4\rho} a^i_l \gs 2^{-1}C^{-1}_4\bar{C}^{-1}_7 p  
$$
and 
$$
\sum_{l\in I(t^j)\backslash I(t^i)} t^j_l1_{k_l> 4\rho} y_l\gs 2^{-1}\bar{C}^{-1}_7\sum_{l\in I(t^j)\backslash I(t^i)} k_l1_{k_l> 4\rho} a^i_l \gs 2^{-1}C^{-1}_4\bar{C}^{-1}_7 p.
$$
Hence if $\bar{C}_2\gs 2 C_4 \bar{C}_7$ we have contradiction with (\ref{onion}).
\smallskip

\noindent
The above construction is valid till $|K|\gs \left(\frac{ckpe}{v}\right)^v$. Since we know that 
$|K|\gs e^{-p}N\gs \exp(-p+\bar{C}_0 p\log(1+p))$ and $N^{\frac{1}{v}}\gs \exp(-\bar{C}_0\log(1+p))=p^{\bar{C}_{0}}$
it implies that we can choose
$$
k=[e^{-1}\frac{v}{cp}e^{-\frac{p}{v}} N^{\frac{1}{v}}]\gs p^{\bar{C}_{5}}
$$ 
for a large enough $\bar{C}_{5}$. It proves that we can find a set $L_k$ such that $|L_k|=k\gs p^{\bar{C}_{5}}$
where $\bar{C}_{5}$ is large enough and for each $i\in L_k$ there exists $l(i)$ for which $y_{l(i)}\ls \bar{C}_7^{-1}a^i_{l(i)}$.
Moreover $l(i)\neq l(j)$ for $i\neq j$, $i,j\in L_k$. Consequently using that $a^i_l\gs C_5^{-1}$ we get
$$
\|z\|^2_2\gs \sum_{i\in L_k}(x_i-y_i)^2\gs (\bar{C}^{-1}_3-\bar{C}^{-1}_7)^2(\sum_{i\in L_k} a^i_{l(i)})^2\gs \bar{C}_8^{-2}|k|, 
$$
where $\bar{C}_8=(\bar{C}^{-1}_3-\bar{C}^{-1}_7)^{-1}C_5$. Therefore 
$$
\|z\|_2=\|x-y\|_2\gs \bar{C}_8^{-1} k^{\frac{1}{2}}\gs \bar{C}_8^{-1}p^{2^{-1}\bar{C}_5}.
$$
Hence by $\|z\|_2\gs p^{\bar{C}_4}$ with $\bar{C}_{4}$ sufficiently large.
It proves (\ref{knot}).
\smallskip

\noindent 
As we have mentioned we can bound $\beta$ in (\ref{n6}) by $\bar{C}_1 p$.
Therefore (\ref{n6}) implies that
$$
\P(M>n_0)=\P(\exists K, |K|>n_0:\;Y\in S_K)\ls e^{-u},\;\;\mbox{for}\;u=\beta^{-1}p^{\bar{C}_4}\gs p^{\bar{C}_4-\bar{C}_1}.
$$
Since $\bar{C}_4$ can be sufficiently large we may require that $\bar{C}_4-\bar{C}_1>1$ and hence we 
have a contradiction
$$
e^{-p}\ls \P(M>n_0)<e^{-p}.
$$
This proves the result.
\end{dwd}

\end{document}